\documentclass[10pt]{amsart}
\usepackage{amsmath,amstext,amssymb,amsopn,amsthm,mathrsfs}
\usepackage{graphicx}
\usepackage[english]{babel}
\usepackage{times}
\usepackage{fancybox}
\usepackage{epsfig}
\usepackage{graphics}
\usepackage{color}
\usepackage{amsmath}
\usepackage{amsfonts}
\usepackage[ansinew]{inputenc}
\usepackage{tikz}
\usepackage{fancyhdr}
\usepackage{lastpage}
\usepackage{enumerate}
\usepackage{graphicx}
\usepackage{times,multicol,amsmath}
\usepackage{array}
\usepackage{fancyhdr}
\usepackage{lastpage}
\usepackage{enumerate}
\usepackage{graphicx}
\usepackage{amsmath}
\usepackage{amsfonts}
\usepackage{amssymb}
\newtheorem{theorem}{Theorem}[section]

\numberwithin{equation}{section}
\newcommand{\dem}{\medskip \par \noindent \mbox{\bf Proof. }}
\def\ep{\hfill{$\Box $}}

\begin{document}

\title[On the Cantor set and the Cantor-Lebesgue function] {On the Cantor set and the Cantor-Lebesgue functions}

\author{Lihang Liu}
\address{Department of Mathematics, University of New Mexico, Albuquerque, NM.}
\email{[Lihang Liu]leonabouttime@unm.edu}
\author{Wilfredo O. Urbina}
\email{[Wilfredo Urbina]wurbinaromero@roosevelt.edu}

\begin{abstract}
The ternary Cantor set $\mathcal{C}$, constructed by George Cantor in 1883, is the best known example of a perfect nowhere-dense set in the real line. The present article we study the basic properties  $\mathcal{C}$ and also study in detail the  ternary expansion characterization   $\mathcal{C}$. We then consider the Cantor-Lebesgue function defined on $\mathcal{C},$ prove its basic properties and study its continuous extension to $[0,1].$ We also consider the  geometric construction of $F$ as the uniform limit of polygonal functions. Finally, we consider the Lebesgue's function defined  from $\mathcal{C}$ onto $[0,1]^2$ and onto $[0,1]^3,$ as well as their continuous extension to $[0,1],$ i.e., obtained the Lebesgue's space filling curves. Finally we discuss Hausdorff's theorem, which is a natural  generalization of the definition of Lebesgue's functions, that states that any compact metric space is a continuous image of the Cantor set $\mathcal{C}$. This notes are the outgrowth of a (zoom)-seminar during the 2020 spring and summer semesters based on H. Sagan's book \cite{sagan}, thus there is not any pretension of originality and/or innovation in this notes and its main objective is a systematic exposition of these topics, in other words its objective is being a review.
\end{abstract}
\maketitle

\section{Cantor Ternary Set.}
The Cantor  ternary set $\mathcal{C}$ is probably the best known example of a {\em perfect nowhere-dense} set in the real line. It was constructed by George Cantor in 1883, see \cite{cantor}, but it was not the first perfect nowhere-dense set  in the real line to be constructed. The first construction was done by the a British mathematician Henry J. S. Smith in 1875, but not many mathematicians were aware of Smith's construction. Vito Volterra, still a graduate student in Italy, also showed how to construct such a set  in 1881 (his famous Volterra set\footnote{In \cite{dimartinourb1} the Volterra set is studied in detail.}$\mathcal{V}$,that allowed him  to construct a counter-example of a function with bounded derivative that exists everywhere but the derivative is not Riemann integrable in any closed bounded interval i.e., the Fundamental Theorem of Calculus fails, see \cite[chapter 4]{bres}) but he published his result in an Italian journal not widely read. In 1883 Cantor rediscover this construction himself, and due to Cantor's prestige, the Cantor ternary set has become the typical example of a perfect nowhere-dense set. Following D. Bresoud \cite{bres}, we will refer as Smith-Volterra-Cantor sets, denoted as $SVC(m)$ sets, to the family of examples of perfect, nowhere-dense sets exemplified by the work of Smith, Voterra and Cantor, obtained by removing in the $n$-iteration, an open interval of length $1/m^n$ from the center of the remaining closed intervals. Observe that $\mathcal{C} = SVC(3)$.\\

As it is well known, $\mathcal{C}$ is obtained from the closed interval $[0,1]$ by a sequence of deletions of open intervals known as ``middle thirds".
We begin with the interval $[0,1]$, let us call it $\mathcal{C} _0$, and remove the middle third, leaving us with the union of two closed intervals of length $1/3$,
$$\mathcal{C}_1= \left[0,\frac{1}{3}\right]  \cup \left[\frac{2}{3}, 1\right] .$$ 
Now we remove the middle third from each of these intervals, leaving us with the union of four closed intervals of length $1/9$,
\begin{equation*}
\mathcal{C}_2=\left[0,\frac{1}{9}\right] \cup \left[\frac{2}{9},\frac{1}{3}\right] \cup\left[\frac{2}{3},\frac{7}{9}\right] \cup\left[\frac{8}{9},1\right] .
\end{equation*}
Then we remove the middle third of each of these intervals leaving us with eight intervals of length $1/27$,
\begin{equation*}
\mathcal{C}_{3}= \left[0,\frac{1}{27}\right] \cup\left[\frac{2}{27}, \frac{1}{9}\right] \cup\left[\frac{2}{9}, \frac{7}{27}\right] \cup \left[ \frac{8}{27},\frac{1}{3}\right]  \cup \left[\frac{2}{3},\frac{19}{27}\right] \cup\left[\frac{20}{27},\frac{7}{9}\right] \cup\left[\frac{8}{9},\frac{25}{27}\right] \cup\left[\frac{26}{27},1\right].
\end{equation*}

We continue this process inductively, then for each $n=1,2, 3\cdots $ we get a set $\mathcal{C} _n$  which is the union of $2^n$ closed intervals of length $1/3^n$. This iterative construction is illustrated in the following figure, for the first four steps:
\begin{figure}[!h]
\begin{center}
\includegraphics[width=3in]{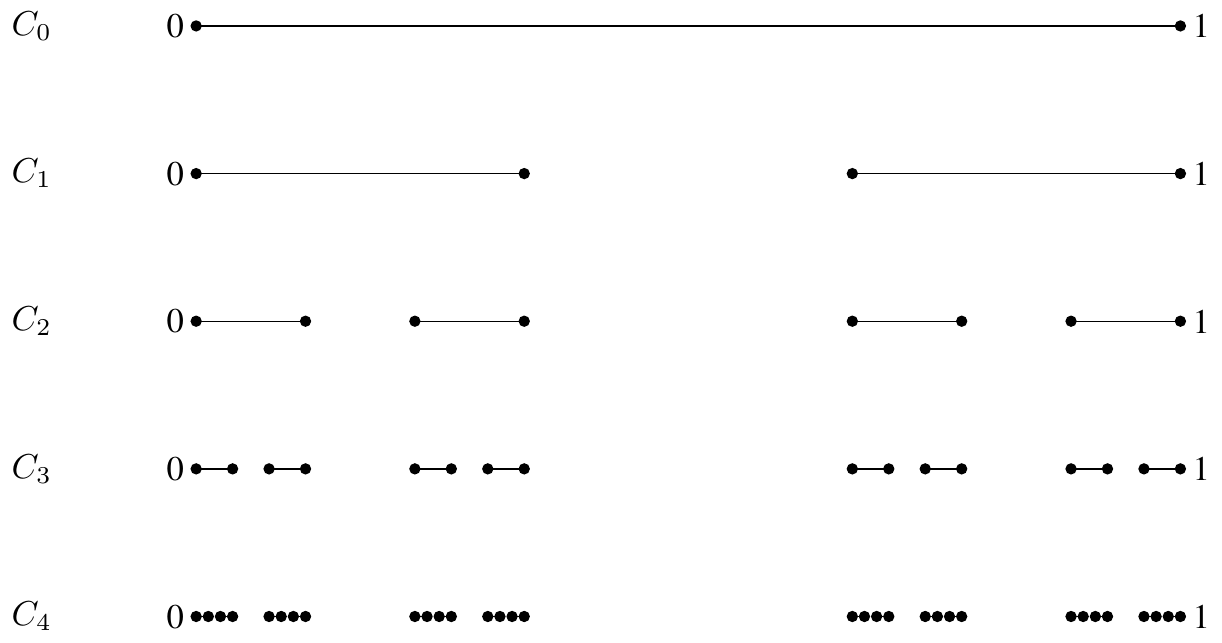}
\end{center}
\caption{The Cantor set $\mathcal{C}$ first four steps}
\end{figure}

Then, the {\em Cantor ternary set} $\mathcal{C}$ is defined as the intersection
\begin{equation}
\mathcal{C}= \bigcap_{n=0}^\infty \mathcal{C}_n.
\end{equation}

\begin{itemize}
\item  Clearly $\mathcal{C} \neq \emptyset $, since  trivially $0,1 \in \mathcal{C}$ and moreover,  if $y$ is the end point of some closed subinterval of a given  $\mathcal{C}_n$ then it is also the end point of some of the subintervals of $\mathcal{C}_{n+1}$. Because at each stage, endpoints are never removed, it follows that $y \in \mathcal{C}_n$ for all $n$. Thus $\mathcal{C}$ contains all the end points of all the intervals that make up each of the sets $\mathcal{C}_n$ (alternatively,  the endpoints to the intervals removed)  all of which are rational ternary numbers in $[0,1]$, i.e., rational triadic, numbers of the form $k/3^n$. 
\item Moreover,  $\mathcal{C}$ is a a closed set, being the countable intersection of closed sets, and  trivially bounded, since it is a subset of $[0,1]$. Therefore, by the Heine-Borel theorem $\mathcal{C}$ is a {\em compact set}.
\item Additionally, $\mathcal{C}$ is a {\em perfect set}, since it is, as we already know closed, and  that every point of $\mathcal{C}$ is a limit point i.e., it is approachable arbitrarily closely by the endpoints of the intervals removed (thus for any $x\in \mathcal{C}$ and for each $n \in \mathbb{N}$ there is an endpoint, let us call it $y_n \in \mathcal{C}_n$, such that $|x-y_n| < 1/3^n$).
\item $\mathcal{C}$ is a {\em nowhere-dense} set, that is, there are no intervals included in $\mathcal{C}$. One way to prove that is taking two arbitrary points in $\mathcal{C}$ we can always find a number between them that requires the digit 1 in its ternary representation, and therefore there are no intervals included in $\mathcal{C}$, thus $\mathcal{C}$ is a nowhere dense set. Alternatively, we can prove this simply by contradiction. Assuming that there is an interval $I=[a,b] \subset \mathcal{C}, \; a <b$. Then $I=[a,b] \subset \mathcal{C}_n$ for all $n$ but as $|\mathcal{C}_n| \rightarrow 0$ as $n \rightarrow \infty$ then $|I| = b-a =0.$
\item Finally, $\mathcal{C}$ has measure zero, since its Lebesgue measure can be obtained  calculating the measure of its complement in $[0,1]$ which is the sum of the length of all  open intervals removed in constructing it, thus
$$ m(\mathcal{C}) = 1 - \sum_{n=1}^\infty \frac{2^{n-1}}{3^n} = 1 -\frac{1}{3}  \sum_{n=0}^\infty (\frac{2}{3})^{n} = 1 -  \frac{1/3}{1- 2/3}=1-1=0.$$\\

\end{itemize}

This construction of $\mathcal{C}$, removing a fixed proportion of each subinterval in each of the iterative steps, give us another important property of the Cantor set, its self-similarity (fractal characteristic) across scales, this is illustrated in the following figure:

\begin{figure}[!h]
\begin{center}
\includegraphics[width=3in]{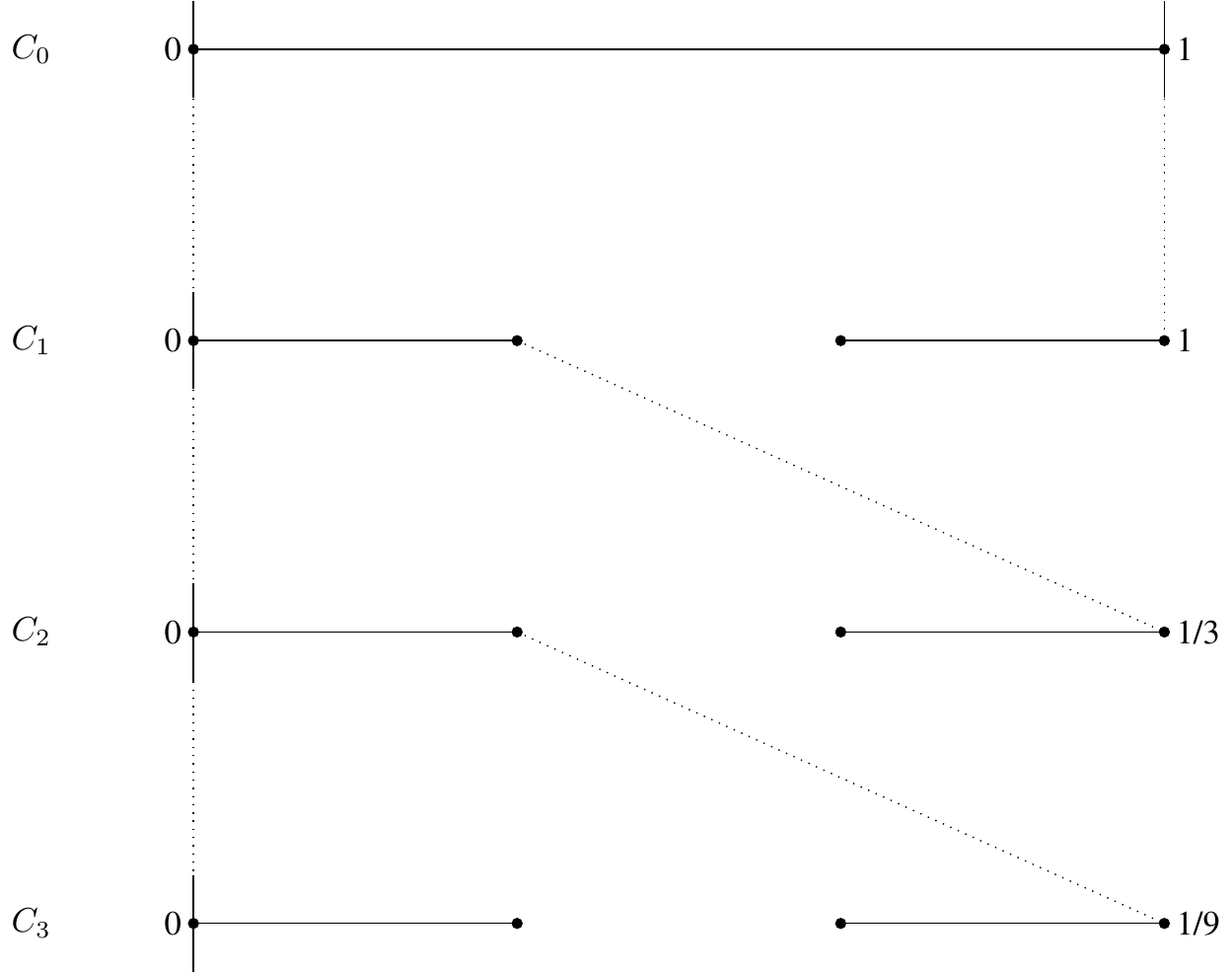}
\end{center}
\caption{The Cantor set $\mathcal{C}$ proportionality}
\end{figure}
In order to compute the Hausdorff dimension of $\mathcal{C}$, S. Abbott in his book used a nice small trick (avoiding the technicalities of the definition of Hausdorff dimension) by dilating $\mathcal{C}_1$ by a the factor $3$ (for more details of this argument see \cite{abot} page 77) obtaining
$$ [0,1] \cup [2, 3].$$
Thus if we continue the iterations we will get two Cantor-like sets, then it can be concluded that
$$3^d= 2, $$
and therefore the Hausdorff dimension of $\mathcal{C}$ is
 $$d= \frac{\ln 2}{\ln 3}=0.630930.$$

There is an alternative characterization of $\mathcal{C}$, the {\em ternary expansion characterization}.
Consider the ternary representation for  $x \in [0,1],$
 
\begin{equation}\label{ternaryexp}
 x = \sum_{k=1}^\infty \frac{\varepsilon_k(x)}{3^k}, \quad \varepsilon_k(x)=0, 1, 2 \quad \mbox{for all} \,k = 1, 2, \cdots 
\end{equation}

Observe that removing the elements where at least one of the \(\varepsilon_{k}\) is equal to one is the same as removing the middle third in the iterative construction, thus 
the Cantor ternary set is the set of numbers in $[0,1]$ that can be written in base 3 without using the digit 1, i.e.,
\begin{equation}\label{ternaryChar}
\mathcal{C}=\left\{x\in [0,1] : x = \sum_{k=1}^\infty \frac{\varepsilon_k(x)}{3^k}, \quad \varepsilon_k(x)=0, 2 \quad \mbox{for all} \; k = 1, 2, \cdots
 \right\}.
\end{equation}

 We want to prove this formally. First we will briefly discuss the binary and ternary representation of real numbers. It will suffice to consider real numbers between 0 and 1, since the representations for other real numbers can then be obtained by adding a positive or negative number.

If $x \in[0,1],$ we will use a repeated bisection procedure to associate a sequence $\left(\eta_{n}\right)$ of 0's and  1's as follows. If $x \neq \frac{1}{2}$ belongs to the left subinterval $\left[0, \frac{1}{2}\right]$ we take $\eta_{1}:=0,$ while if $x$ belongs to the right subinterval $\left[\frac{1}{2}, 1\right]$ we take $\eta_{1}:=1.$ If $x=\frac{1}{2},$ then we may take $\eta_{1}$ to be either 0 or $1 .$ In any case, we have
\[
\frac{\eta_{1}}{2} \leq x \leq \frac{\eta_{1}+1}{2}.
\]
We now bisect the interval $\left[\frac{1}{2} \eta_{1}, \frac{1}{2}\left(\eta_{1}+1\right)\right].$ If $x$ is not the bisection point and belongs to the left subinterval we take $\eta_{2}:=0,$ and if $x$ belongs to the right subinterval we take $\eta_{2}:=1 .$ If $x=\frac{1}{4}$ or $x=\frac{3}{4},$ we can take $\eta_{2}$ to be either $0$ or $1.$ In any case, we have
\[
\frac{\eta_{1}}{2}+\frac{\eta_{2}}{2^{2}} \leq x \leq \frac{\eta_{1}}{2}+\frac{\eta_{2}+1}{2^{2}}.
\]
We continue this bisection procedure, assigning at the $n$-th stage the value $\eta_{n}:=0$ if $x$ is not the bisection point and lies in the left subinterval, and assigning the value $\eta_{n}:=1$ if $x$ lies in the right subinterval. In this way we obtain a sequence $\left(\eta_{n}\right)$ of 0's or 1's that correspond to a nested sequence of intervals containing the point $x .$ For each $n,$ we have the inequality
\begin{equation}\label{ineq1}
\frac{\eta_{1}}{2}+\frac{\eta_{2}}{2^{2}}+\dots+\frac{\eta_{n}}{2^{n}} \leq x \leq \frac{\eta_{1}}{2}+\frac{\eta_{2}}{2^{2}}+\cdots+\frac{\eta_{n}+1}{2^{n}}.
\end{equation}
If $x$ is the bisection point at the $n$-th stage, then $x=m / 2^{n}$ with $m$ odd. In this case, we may choose either the left or the right subinterval; however, once this subinterval is chosen, then all subsequent subintervals by the bisection procedure are determined. For instance, if we choose the left subinterval so that $\eta_{n}=0,$ then $x$ is the right endpoint of all subsequent subintervals, and hence $\eta_{k}=1$ for all $k \geq n+1 .$ On the other hand, if we choose the right subinterval so that $\eta_{n}=1,$ then $x$ is the left endpoint of all subsequent subintervals, and hence $\eta_{k}=0$ for all $k \geq n+1 .$ For example, if $x=\frac{3}{4},$ then the two possible expansions for $x$ are $0._210111, \ldots$ and $0._21100, \ldots .$

To summarize: If $x \in[0,1],$ then there exists a sequence $\left(\eta_{n}\right)$ of 0's and 1's such that inequality above holds for all $n \in \mathbb{N} .$ In this case we write
\[
x= 0._2  \eta_{1} \eta_{2} \cdots \eta_{n} \cdots
\]
and call it the binary representation of $x .$ This representation is unique except when $x=m / 2^{n}$ for $m$ odd, in which case $x$ has the two representations
\[
x=0._2\eta_{1} \eta_{2} \cdots \eta_{n-1} 1000 \cdots=0._2 \eta_{1} \eta_{2} \cdots \eta_{n-1} 0111 \cdots
\]
one ending in 0's and the other ending in 1's. Conversely, each sequence of 0's and 1's is the binary representation of a unique real
number in $[0,1] .$ The inequality (\ref{ineq1}) determines a closed interval with length $1 / 2^{n}$ and the sequence of these intervals is nested. Therefore, the nested interval property implies that there exists a unique real number $x$ satisfying the inequality for every $n \in \mathbb{N}$. Consequently, $x$ has the binary representation 
$$x= 0._2  \eta_{1} \eta_{2} \cdots \eta_{n} \cdots\\$$

Analogously, if $x \in[0,1],$ we will use a repeated trisection procedure to associate a sequence $\left(\varepsilon_{n}\right)$ of 0's, 1's and 2's as follows. If $x \neq \frac{1}{3}$ belongs to the left subinterval $\left[0, \frac{1}{3}\right]$ we take $\varepsilon_{1}:=0,$ if $x$ belongs to the right subinterval $\left[\frac{1}{3}, \frac23\right]$ we take $\varepsilon_{1}:=1$ while if $x$ belongs to the right subinterval $\left[\frac{2}{3}, 1\right]$ we take $\varepsilon_{1}:=2 .$ If $x=\frac{1}{3},$ then we may take $\varepsilon_{1}$ to be either 0 or $1 ,$ and if   $x=\frac{2}{3},$ then we may take $\varepsilon_{1}$ to be either 1 or 2. In any case, we have
\[
\frac{\varepsilon_{1}}{3} \leq x \leq \frac{\varepsilon_{1}+1}{3}
\]
We now trisect the interval $\left[\frac{1}{3} \varepsilon_{1}, \frac{1}{3}\left(\varepsilon_{1}+1\right)\right] .$ If $x$ is not the trisection point and belongs to the left subinterval we take $\varepsilon_{2}:=0,$ and if $x$ belongs to the right subinterval we take $\varepsilon_{2}:=1 .$ If $x=\frac{1}{9}, \frac29, \frac49, \frac59, \frac79$ or $\frac{8}{9},$ we can take $\varepsilon_{2}$ to be either 0, 1 or $2 .$ In any case, we have
\[
\frac{\varepsilon_{1}}{3}+\frac{\varepsilon_{2}}{3^{2}} \leq x \leq \frac{\varepsilon_{1}}{3}+\frac{\varepsilon_{2}+1}{3^{2}}
\]
We continue this trisection procedure, assigning at the $n$ th stage the value $\varepsilon_{n}:=0$ if $x$ is not the trisection point and lies in the left subinterval, assigning the value $\varepsilon_{n}:=1$ if $x$ lies in the middle subinterval and assigning the value $\varepsilon_{n}:=2$ if $x$ lies in the right subinterval. In this way we obtain a sequence $\left(\varepsilon_{n}\right)$ of 0's, 1's or 2's that correspond to a nested sequence of intervals containing the point $x .$ For each $n,$ we have the inequality
\begin{equation}\label{ineq2}
\frac{\varepsilon_{1}}{3}+\frac{\varepsilon_{2}}{3^{2}}+\dots+\frac{\varepsilon_{n}}{3^{n}} \leq x \leq \frac{\varepsilon_{1}}{3}+\frac{\varepsilon_{2}}{3^{2}}+\cdots+\frac{\varepsilon_{n}+1}{3^{n}}.
\end{equation}
If $x$ is the trisection point at the $n$-th stage, then $x=m / 3^{n}$. In this case, we may choose either the left or the right subinterval; however, once this subinterval is chosen, then all subsequent subintervals in the trisection procedure are determined. [For instance, if we choose the left subinterval so that $\varepsilon_{n}=0,$ then $x$ is the right endpoint of all subsequent subintervals, and hence $\varepsilon_{k}=1$ for all $k \geq n+1 .$ On the other hand, if we choose the right subinterval so that $\varepsilon_{n}=2,$ then $x$ is the left endpoint of all subsequent subintervals, and hence $\varepsilon_{k}=0$ for all $k \geq n+1 .$ For example, if $x=\frac{2}{3},$ then the two possible sequences for $x$ are $1,2,2,2,2, \ldots$ and $ 2,0,0,0, \ldots .$

To summarize: If $x \in[0,1],$ then there exists a sequence $\left(\varepsilon_{n}\right)$ of 0's, 1's and 2's such that inequality above holds for all $n \in \mathbb{N} .$ In this case we write
\[
x= 0._3  \varepsilon_{1} \varepsilon_{2} \cdots \varepsilon_{n} \cdots
\]
and call it the {\em ternary representation} of $x .$ This representation is unique except when $x=m / 3^{n}$ for  some $m,$ in which case $x$ has the two representations
\[
x=0._{3} \varepsilon_{1} \varepsilon_{2} \cdots \varepsilon_{n-1} 1000 \cdots=0_{3}. \varepsilon_{1} \varepsilon_{2} \cdots \varepsilon_{n-1} 0222 \cdots
\]
one ending in 0's and the other ending in 2's. Conversely, each sequence of 0s and 1s is the binary representation of a unique real
number in $[0,1] .$ The inequality above determines a closed interval with length $1 / 3^{n}$ and the sequence of these intervals is nested. Therefore, the nested interval property implies that there exists a unique real number $x$ satisfying the inequality for every $n \in \mathbb{N}$. Consequently, $x$ has the ternary representation 
$$x= 0._3  \varepsilon_{1} \varepsilon_{2} \cdots \varepsilon_{n} \cdots$$
Thus, for the dyadic rational numbers $m/2^n$ there are two possible dyadic expansions as
$$ \frac{m}{2^n} = \frac{m-1}{2^n} + \frac{1}{2^n}  =  \frac{m-1}{2^n} + \sum_{k=n+1}^\infty \frac{1}{2^k}.$$ 
Similarly,  for the ternary rational  numbers $m/3^n$ there are two possible ternary expansions, since
$$ \frac{m}{3^n} = \frac{m-1}{3^n} + \frac{1}{3^n}  =  \frac{m-1}{3^n} + \sum_{k=n+1}^\infty \frac{2}{3^k}.$$ 
Therefore the dyadic and the ternary representations  are unique except for the dyadic  or  ternary rational numbers, in that case we will take  the infinite expansions representations for the dyadic and ternary rational numbers.\\

Now, let us prove that  $x \in \mathcal{C}$ if and only if  $x=\sum_{k=1}^{\infty} \frac{\varepsilon_{k}}{3^{k}},\; \text { where } \varepsilon_{k}=0,\text { or } 2.$\\
\begin{itemize}
\item  We have seen that, for $n=1,$ an interval $\left(a_{n}, b_{n}\right)$ is removed if $a_{1}=0._{3} 1=0._{3} 0\overline{2}$ and $b_{1}=0._{3} 2 .$ For $n=2,$ it is removed if $a_{2}=0._{3} 01= 0._300\overline{2},\; b_{2}=0._{3} 02$ or $a_{2}=0._{3} 21=0._320\overline{2},\; b_{2}=0._{3} 22 .$ In general, it is removed if and only if $a_{n}, b_{n}$ can be written in the form
\begin{eqnarray*}
a_{n}&=&\sum_{k=1}^{n-1} \frac{\varepsilon_{k}}{3^k} + \frac1{3^n}= 0._{3}\varepsilon_{1}\varepsilon_{2}\varepsilon_{3} \ldots \varepsilon_{n-1} 1 \\&=&0._{3}\varepsilon_{1}\varepsilon_{2}\varepsilon_{3} \ldots \varepsilon_{n-1} 0\overline{2}= \sum_{k=1}^{n-1} \frac{\varepsilon_{k}}{3^k} + \sum_{k=n+1}^{\infty} \frac{2}{3^k},\\
b_{n}&=&\sum_{k=1}^{n-1} \frac{\varepsilon_{k}}{3^k} + \frac2{3^n} = 0._{3}\varepsilon_{1}\varepsilon_{2}\varepsilon_{3} \ldots\varepsilon_{n-1} 2.
\end{eqnarray*}

Suppose this is true for a given $n$, then for the $(n+1)$-th step, we have
\begin{eqnarray*}
a_{n+1} &=&b_{n}+\frac1{3^{n+1}}= 0._{3}\varepsilon_{1}\varepsilon_{2}\varepsilon_{3} \ldots \varepsilon_{n-1} 2 +\frac1{3^{n+1}} \\
&=&0._{3}\varepsilon_{1}\varepsilon_{2}\varepsilon_{3} \ldots \varepsilon_{n-1}  21= 0._{3}\varepsilon_{1}\varepsilon_{2}\varepsilon_{3} \ldots \varepsilon_{n-1}  20\overline{2}\\
b_{n+1} &=&b_n+\frac2{3^{n+1}} = 0._{3}\varepsilon_{1}\varepsilon_{2}\varepsilon_{3} \ldots\varepsilon_{n-1} 22
\end{eqnarray*}
or 
\begin{eqnarray*}
a_{n+1} &=&a_{n}-\frac{2}{3^{n+1}}=0._{3}\varepsilon_{1}\varepsilon_{2}\varepsilon_{3} \ldots\varepsilon_{n-1} 1-\frac{2}{3^{n+1}} \\
&=&0._{3}\varepsilon_{1}\varepsilon_{2}\varepsilon_{3} \ldots\varepsilon_{n-1} 01 = 0._{3}\varepsilon_{1}\varepsilon_{2}\varepsilon_{3} \ldots \varepsilon_{n-1} 00\overline{2}\\
b_{n+1} &=&a_{n}-\frac{1}{3^{n+1}}=0._{3}\varepsilon_{1}\varepsilon_{2}\varepsilon_{3} \ldots\varepsilon_{n-1} 1-\frac{1}{3^{n+1}} \\
&=&0._{3}\varepsilon_{1}\varepsilon_{2}\varepsilon_{3} \ldots\varepsilon_{n-1} 02
\end{eqnarray*}

Therefore by the induction principle we conclude that  
$$x=\sum_{k=1}^{\infty} \frac{\varepsilon_{k}}{3^{k}}, \quad \text { where } \varepsilon_{k}=0\text { or } 2.$$
Observe that this implies, in particular, that the left end points of the disjoint intervals constituting $\mathcal{C}_{n}$ have a finite triadic representation with ending digit $2$; more precisely, they can be written as 
$\sum_{k=1}^{\infty} \frac{\varepsilon_{k}}{ 3^{k}}$ where $\varepsilon_{n} \in\{0,2\}$ for $0<k \leq n,\; \varepsilon_{k} =2$ and $\varepsilon_{k}=0$ for $k>n,$ (this can also be proved directly using mathematical induction) and the right end points of the disjoint intervals constituting $\mathcal{C}_{n}$ have a infinite periodic triadic representation with period $2$. Moreover, observe that $\mathcal{C}_{n}$ consists of exactly $2^{n}$ disjoint intervals, which is exactly the number of points in [0,1] of the form $\sum_{k=1}^{\infty} \frac{\varepsilon_{k}}{3^{k}},$ where $a_{k} \in\{0,2\}$ for $0<k \leq n$ and $a_{k}=0$ for $k>n .$\\
\item  Now, if  $x \in\left(a_{n}, b_{n}\right),$ then
$$x=0._{3}\varepsilon_{1}\varepsilon_{2} \varepsilon_{3} \ldots\varepsilon_{n-1} 1 \tau_{n+1} \tau_{n+2} \ldots $$
If all $\tau_{n+j}=0,$ then necessarily $x=a_{n} \notin\left(a_{n}, b_{n}\right),$ or if all $\tau_{n+j}=2,$ then  necessarily $x=b_{n} \notin$ $\left(a_{n}, b_{n}\right).$ Therefore, at least one of the $\tau_{n+j} \neq 0$ and at least one $\tau_{n+j} \neq 2.$
 
Hence, a point $x$ is removed in the course of the construction of $\mathcal{C}$ if and only if both its ternary representations (if there are two) contain at least one digit $1.$\\

\end{itemize}

For a systematic treatment  of the Cantor set and Cantor like sets we refer to \cite{dimartinourb}.

\section{The Cantor-Lebesgue function}

The {\em Cantor-Lebesgue} function, also called  simply {\em Cantor function}, is defined on $\mathcal{C}$ by
\begin{equation}
F(x)=\sum_{k=1}^{\infty} \frac{\eta_{k}}{2^{k}} \quad \text { if }\; x\in \mathcal{C} \; \text { with }\; x=\sum_{k=1}^{\infty}  \frac{\varepsilon_{k}}{3^{k}}, \text { where }\eta_{k}=\varepsilon_{k} / 2,
\end{equation}
where we need to choose the expansion of $x$ in which $\varepsilon_{k}=0$ or $2.$ 
Observe that the action of the function $F$ can be written as 
\[
F\left(0._{\dot{3}}\varepsilon_{1}\varepsilon_{2}\varepsilon_{3}\varepsilon_{4} \ldots\right)=0._{\dot{2}} \eta_{1} \eta_{2} \eta_{3} \eta_{4} \ldots, \eta_{k}=0 \, \text { or } 1.
\]
\begin{itemize}
\item $F$ is well defined since from the discussion above, for almost all points of $[0,1]$ (except for the binary and ternary rationals, $m / 2^{n},\; m / 3^{n}$ ) have a unique representation. For the case of the  binary and triadic rationals they have two representations, in order to apply the definition of $F(x),$ we need the infinite representation with period $2$
$$x=\frac{m}{ 3^{n}} =  0._{3} \varepsilon_{1} \varepsilon_{2} \cdots \varepsilon_{n-1} 1000 \cdots=0_{3}. \varepsilon_{1} \varepsilon_{2} \cdots \varepsilon_{n-1} 0222 \cdots,$$
\begin{eqnarray*}
F(x) &=& 0_{2}. \left(\varepsilon_{1}/2\right)\left( \varepsilon_{2}/2\right)\cdots \left( \varepsilon_{n-1}/2\right) 0111 \cdots \\
&=& \left(\varepsilon_{1}/2\right)\left( \varepsilon_{2}/2\right) \cdots \left(\varepsilon_{n-1}/2\right) 1000   \cdots.
\end{eqnarray*}
\item $F: \mathcal{C} \rightarrow[0,1]$ is not injective, simply observe that both $0._{3} 0\overline{2}$ and $0._{3} 2,$ are mapped into $0._{2} 1=0._{2} 0\overline{1}$. But $F$ is surjective, as given $y \in  [0,1]$ with  binary expansion of the form $\sum_{k=1}^\infty \frac{\eta_k}{2^k}$ with $\eta_k \in \{0,1\},$ then considering $x =\sum_{k=1}^\infty \frac{2 \eta_k}{3^k}$ 
it is clear that $x \in \mathcal{C}$ and such that by construction $F(x) =y.$ This proves the fact that $F$ is surjective. 
\item Let us now prove the continuity of $F$ on $\mathcal{C}$. Suppose that $x_0 \in  \mathcal{C},$ fix, and let  $x \in  \mathcal{C},$ such that $\left|x-x_{0}\right|<\frac{1}{3^{2 n}} .$ Then, $x_{0}$ and $x$ cannot differ in the first $2 n$ ternary places and we have
\[
\begin{array}{c}
x_{0}=0._{3}\varepsilon_{1}\varepsilon_{2} \ldots\varepsilon_{2 n-1}\varepsilon_{2 n}\varepsilon_{2 n+1} \ldots \\
x=0._{3}\varepsilon_{1}\varepsilon_{2}\ldots\varepsilon_{2 n-1}\tau_{2 n}\tau_{2 n+1} \ldots
\end{array}
\]
Suppose, to the contrary, that $\varepsilon_{2 n} \neq \tau_{2 n} .$ Then, $\left|\varepsilon_{2 n}-\tau_{2 n}\right|=2$ and
\[
x-x_{0}= (\tau_{2 n}-\varepsilon_{2 n} )/ 3^{2 n}+ (\tau_{2 n+1}-\varepsilon_{2 n+1}) / 3^{2 n+1}+\cdots
\]
and, hence,
\[
\begin{aligned}
\left|x-x_{0}\right| & \geq\left(2 / 3^{2 n}\right)-\left(2 / 3^{2 n+1}\right)(1+1 / 3+1 / 9+\cdots) \\
&=2 / 3^{2 n}-1 / 3^{2 n}=1 / 3^{2 n},
\end{aligned}
\]
instead of $<1 / 3^{2 n} .$ \\

Therefore, given $\varepsilon >0$ there exists $N$ such that $1/2^{2N} < \varepsilon,$  and taking $\delta = 1/3^{2N}$ we have
that if $\left|x-x_{0}\right|<\frac{1}{3^{2 N}} $ then
\begin{eqnarray*}
|F(x)-F\left(x_{0}\right)|&=& | \frac{(\tau_{2 N+1}/2-\varepsilon_{2 N+1}/2)}{ 2^{2N+1}}+\frac{(\tau_{2 N+2}/2-\varepsilon_{2 N+2}/2)}{2^{2N+2}}+\cdots|\\
&=& | \frac{(\tau_{2 N+1}-\varepsilon_{2 N+1})}{ 2^{2N+2}}+\frac{(\tau_{2 N+2}-\varepsilon_{2 N+2})}{2^{2N+3}}+\cdots| \\
&\leq& \frac{2}{2^{2N+2}} \Big( 1+\frac12+ \frac14+ \cdots\Big)=\frac1{2^{2N}} < \varepsilon\\
\end{eqnarray*}
hence the continuity of $F$ at $x_0$ follows. 

We have actually  proved that $F$ is uniformly continuous on $\mathcal{C}$, since $\delta$ is independent of $x_0,$ since  given $\epsilon>0$,  taking  $1/2^{2N} < \varepsilon,$ and $\delta = 1/3^{2N}$ we have  for all $x,x' \in \mathcal{C}$ such that $|x-x'|<\delta$ implies  $|F(x)-F(x')| \leq \varepsilon .$ \\
\item Additionally, $F(0)=0$ and $F(1)=1$. Since $0=\displaystyle\sum_{k=1}^{\infty} \frac{0}{ 3^{k}}$ and $1=\displaystyle\sum_{k=1}^{\infty} \frac{2}{ 3^{k}}=0._3\overline{2},$ thus
\[
F(0)=\displaystyle\sum_{k=1}^{\infty} \frac{0}{ 2^{k}}=0 \text { and } F(1)=\displaystyle\sum_{k=1}^{\infty} \frac{1}{2^{k}}= \frac{\frac{1}{2}}{1-\frac{1}{2}}=1.\\
\]
\item The function $F$ can be extended to a continuous function on $[0,1]$ as follows. First, we will prove that if $(a_n, b_n)$ is an open interval of the complement of $\mathcal{C},$ then $F(a_n)=F(b_n).$ Hence we may define $F$ to have the constant value $F(a_n)$ in that interval. Then, the definition of $F$ may be extended into all of  $[0,1]$ by defining it on $\mathcal{C}^{c}$ as follows: If $x \in \mathcal{C}^{c},$ then $x\in(a_n, b_n)$ where the open interval $(a_n, b_n)$ is one of those that has been removed from $[0,1]$ in the construction of the Cantor set. Then, necessarily
\[
\begin{array}{l}
a_{n}=0._{3}\varepsilon_{1}\varepsilon_{2}\varepsilon_{3} \ldots \varepsilon_{n-1} 1=0._{3}\varepsilon_{1}\varepsilon_{2}\varepsilon_{3} \ldots \varepsilon_{n-1} 0 \overline{2} \\
b_{n}=0._{3}\varepsilon_{1}\varepsilon_{2}\varepsilon_{3} \ldots\varepsilon_{n-1} 2,
\end{array}
\]
and therefore,
\begin{eqnarray*}
F(a_{n}) &=& 0._{2}(\varepsilon_{1}/2)(\varepsilon_{2}/2)(\varepsilon_{3}/2)\ldots (\varepsilon_{n-1}/2) 0 \overline{1}\\
& =& 0._{2}(\varepsilon_{1}/2)(\varepsilon_{2}/2)(\varepsilon_{3}/2)\ldots (\varepsilon_{n-1}/2) 1 = F(b_{n}).
\end{eqnarray*}

Now, to prove that the extension, also denoted by $F,$ is continuous on $[0,1]$, i.e., we need to prove that for any $x_0 \in [0,1]$ given $\varepsilon >0$  there exist $\delta >0$  such  that if  $x \in [0,1]$ such that $|x-x_0|<\delta$ then  $|F(x)-F(x_0)| \leq \varepsilon .$ 
\begin{enumerate}
\item[i)] If $x_{0} \in \mathcal{C}^c,$ then $x_0 \in (a_n, b_n)$ where the open interval $(a_n, b_n)$ is one of those that has been removed from $[0,1]$ in the construction of the Cantor set and then as we have prove above $F$ is constant on $(a_n, b_n)$ so  $F$ is trivially at $x_0.$
\item[ii)] If $x_{0} \in \mathcal{C},$ then, as $\mathcal{C}$ is perfect, it is either a left accumulation point, a right accumulation point or an accumulation point. We will only consider the case where $x_{0}$ is a left accumulation point, which makes that it a right endpoint of an interval that has been removed in the construction of $\mathcal{C}.$ From above we know that namely $x_0 = 0._{3}\varepsilon_{1}\varepsilon_{2}\varepsilon_{3} \ldots \varepsilon_{n-1}2.$ Then we know, that  the restriction of $F$ to $\mathcal{C},$ is continuous on it and, since $\mathcal{C}$ is compact, then it is uniformly continuous on $\mathcal{C},$ i.e., for any $\epsilon>0,$ there is $\delta>0$ such that
\[
|F\left(x_{1}\right)-F\left(x_{2}\right)|<\epsilon,
\]
for all $x_{1}, x_{2} \in \mathcal{C}$ for which $\left|x_{1}-x_{2}\right|<\delta.$ \\

Let us now examine $|F\left(x\right)-F\left(x_{0}\right)|.$ There are two possibilities: 
\begin{itemize}
\item $x \in \mathcal{C}$: then, we already know
\[
|F(x)-F\left(x_{0}\right)|<\epsilon \quad \text { for } \quad\left|x-x_{0}\right|<\delta,
\]
\begin{figure}[!h]
\centering
\includegraphics[scale=0.5]{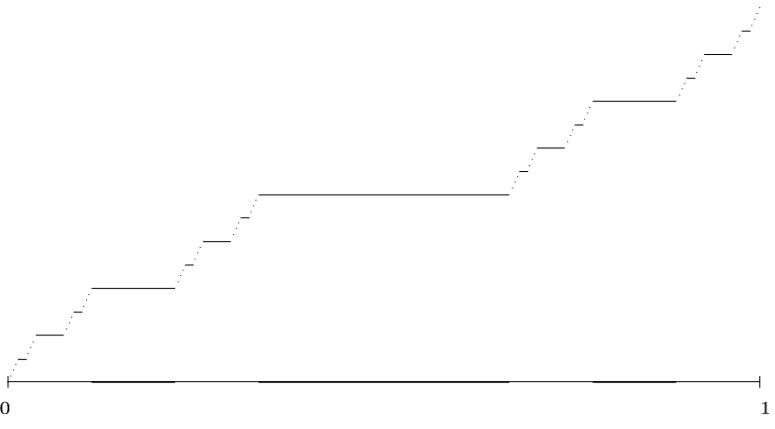}
\caption{The Cantor-Lebesgue function $F$}
\end{figure}
\item  $x \in \mathcal{C}^c,$ then, $x \in(a, b),$ where $(a, b)$ is one of the intervals that have been removed in the construction of $\mathcal{C},$ then  $a, b \in \mathcal{C},\;(a, b) \subset \mathcal{C}^{c} .$ But again, as $F$ is constant on $(a, b),$ then $ F(t) = F(a) = F(b)$ for all $t \in(a, b) .$ Now, let $(a, b)$ be any such interval with $x_{0}<a<b<x_{0}+\delta,$ 
\[
|F(t)-F\left(x_{0}\right)|=|F(a)-F\left(x_{0}\right)| <\epsilon,
\]
for all $t \in(a, b) \subset \mathcal{C}^{c},$ as  $a \in \mathcal{C}.$
\end{itemize}
\end{enumerate}
This, together with the fact that to the left of $x_{0} (x_{0}$ being the right endpoint of a removed interval) which implies the continuity of $F$ at $x_{0} \in \mathcal{C}.$ An analogous proof applies to a right accumulation point of $\mathcal{C},$ and both proofs together take care of a two-sided accumulation point of $\mathcal{C} .$ \\

\item Additionally,  $F$ is nowhere differentiable on $\mathcal{C}$. Let $x, x_n \in \mathcal{C}$ such that
\[
\begin{aligned}
x &=0._3\varepsilon_{1}\varepsilon_{2} \ldots\varepsilon_{ n-1}\varepsilon_{ n}\varepsilon_{ n+1} \varepsilon_{ n+2}\ldots \\
x_{n} &=0._3\varepsilon_{1}\varepsilon_{2} \ldots\varepsilon_{ n-1}\varepsilon_{ n}\tau_{ n+1} \varepsilon_{ n+2} \ldots
\end{aligned}
\]
where $\tau_{2 n+1}=\varepsilon_{2 n+1}+1(\bmod 2) .$ Then
\[
\left|x-x_{n}\right|=2 / 3^{2 n+1}
\]
and then,
\[
F(x)-F\left(x_{n}\right)=\left(t_{2 n+1}-\tau_{2 n+1}\right) / 2^{n+1}
\]
and, hence,
\[
\left|F(x)-F\left(x_{n}\right)\right| /\left|x-x_{n}\right|=\frac{3}{4}\left(\frac{9}{2}\right)^{n} \longrightarrow \infty,
\]
thus, $F$ is nowhere differentiable on $\mathcal{C}.$ Therefore, as $F$ is trivially differentiable on $\mathcal{C}^c,$ with derivative zero, we can conclude that $F$ is a singular function,  i.e., $F'(x)$ is zero a.e., $x \in [0,1]$ without being constant.\footnote{For more on singular functions we refer to \cite{Salem}.} Also it serves as a counterexample to Harnack's extension of the Fundamental Theorem of Calculus to discontinuous functions, which was in vogue at the time (see e.g. \cite[pag. 60]{hawk}).

The graph of the function $F$ is called the {\em devil's staircase}. Moreover, $F$ is an example of an increasing continuous $F$ whose derivative is integrable on $[0,1],$ but  the Fundamental Theorem of Calculus fails: 
$$\int_{0}^{1} F^{\prime}(x) \;dx=0 \neq 1= F(1)-F(0).$$
\item The Cantor-Lebesgue function can also be obtained by a  geometric construction. Let $F_{1}(x)$ be the continuous increasing function on $[0,1]$ defined as
\[
F_{1}(x)=\begin{cases}
(3 / 2) x & \text { for } \;0 \leq x \leq 1 / 3, \\
1 / 2 & \text { for } \;1 / 3<x<2 / 3, \\
(3 / 2) x-1 / 2 & \text { for }\; 2 / 3 \leq x \leq 1.
\end{cases}
\]
 Similarly, let $F_{2}(x)$ be continuous and increasing, defined as
\[
F_{2}(x)=\begin{cases}
(9 / 4) x & \text { if }\; 0 \leq x \leq 1 / 9,\\
1 / 4 & \text { if } \;1 / 9 \leq x \leq 2 / 9, \\
(9 / 4)  x -1/4 & \text { if }\; 2/9 \leq x \leq 1 / 3,\\
1 / 2 & \text { if } 1 / 3\; \leq x \leq 2 / 3,\\
(9 / 4)  x -1 & \text { if } \;2/3 \leq x \leq 7/9,\\
3 / 4 & \text { if }\; 7 / 9 \leq x \leq 8 / 9,\\
(9 / 4)  x -5/4 & \text { if }\; 8/9 \leq x \leq 1.\\
\end{cases}
\]
and $F_{3}(x)$ is also defined to be continuous and increasing, as
\[
F_{3}(x)=\begin{cases}
(27 / 8)  x & \text { if }\; 0 \leq x \leq 1 / 27,\\
1 / 8 & \text { if} \; 1 / 27 \leq x \leq 2 / 27,\\
(27 / 8)  x - 1/8 & \text { if }\; 2 / 27 \leq x \leq 1/9,\\
1 / 4 & \text { if } \; 1 / 9 \leq x \leq 2 / 9,\\
(27 / 8)  x - 1/4 & \text { if }\; 2 / 9 \leq x \leq 7/27,\\
3 / 8 &\text { if } \; 7/ 27 \leq x \leq 8/27,\\
(27 / 8)  x - 5/8 & \text { if }\; 8 / 27 \leq x \leq 1/3,\\
1 / 2 & \text { if } \;1 / 3 \leq x \leq 2 / 3,\\
(27 / 8)  x - 7/4 & \text { if }\; 2 / 3 \leq x \leq 19/27,\\
5 / 8 &\text { if }\; 19 / 27 \leq x \leq 20 / 27,\\
(27 / 8)  x - 15/8 & \text { if }\; 20/27 \leq x \leq 7/9,\\
3 / 4& \text { if  } \;7 / 9 \leq x \leq 8 / 9,\\
(27 / 8)  x -9/4 & \text { if }\; 8/ 9 \leq x \leq 25/27,\\
7 / 8 & \text { if  } \;25 / 27 \leq x \leq 26 / 27, \\
(27 / 8)  x - 19/8 & \text { if }\; 26/27 \leq x \leq 1.
\end{cases}
\]
and so on. \\

Observe $F_{1}$ is continuous, increasing, and constant on the middle third $(1 / 3,2 / 3)=[0,1] \backslash \mathcal{C}_{1}.$  Now observe that $F_{2}$  is linear on $\mathcal{C}_{2}=[0,1 / 9] \cup [2/9,1/3] \cup [2 / 3,7/9]\cup[8/9,1],$ and constant on the intervals  $[1/9,2/9],[1/3.2/3],[7/9,8/9]$ and that can be obtained by imitating this process of flattening out the middle third of each non-constant segment of $F_{1}$. Specifically, it is easy to see that 
\[
F_{2}(x)= \begin{cases}
(1 / 2) F_{1}(3 x) & \text { for } \;0 \leq x \leq 1 / 3 \\
F_{1}(x) & \text { for }\; 1 / 3<x<2 / 3 \\
(1 / 2) F_{1}(3 x-2)+1 / 2 & \text { for } \;2 / 3 \leq x \leq 1
\end{cases}
\]
\begin{figure}[!h]
\centering
\includegraphics[scale=0.95]{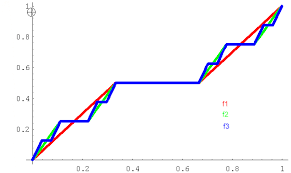}
\caption{The first three approximations $F_1,$ $F_2$, $F_3$ of the Cantor-Lebesgue function $F.$ }
\end{figure}

Iterating this process we observe that $F_{n+1}$ can be obtained by imitating this process of flattening out the middle third of each non-constant segment of $F_{n}$, thus
\[
F_{n+1}(x)= \begin{cases}
(1 / 2) F_{n}(3 x) & \text { for } \;0 \leq x \leq 1 / 3 \\
F_{n}(x) & \text { for }\; 1 / 3<x<2 / 3 \\
(1 / 2) F_{n}(3 x-2)+1 / 2 & \text { for } \;2 / 3 \leq x \leq 1
\end{cases}
\]
and  therefore
 \[
\left|F_{n+1}(x)-F_{n}(x)\right| \leq \frac1{2^{n+1}}.
\]
Moreover, given $m<n$ it is easy to  see that
\[
\left|F_{m}(x)-F_{n}(x)\right| \leq \frac{1}{2^{m}},
\]
as
\begin{eqnarray*}
\left|F_{n}(x)-F_{m}(x)\right| &\leq& \left|F_{m+1}(x)-F_{m}(x)\right| + \cdots +\left|F_{n}(x)-F_{n-1}(x)\right| \\
&\leq& \frac1{2^{m+1}} + \frac1{2^{m+2}}+ \cdots +\frac1{2^{n}}\\
& <&  \frac1{2^{m}} (\frac12+ \frac14+ \cdots + \frac1{2^{n-m}}+ \cdots ) = \frac1{2^{m}}.
\end{eqnarray*}

Observe that this geometric construction allow us to obtain alternative proofs of the properties of $F$:

\item Using the Cauchy criterion for uniform convergence, we conclude that  $\{F_n\}$ is a Cauchy sequence in the uniform norm and therefore $\{F_{n}\}$ converges uniformly and the limit function $F$ is continuous.
\item  $F$ is increasing, since if  $x<y$ be arbitrary points in $[0,1]$ then as $F_{n}(x) \leq F_{n}(y)$ for all $n$ and $F_{n}(x) \rightarrow F(x)$ and $F_{n}(y) \rightarrow F(y),$ using the order limit theorem we conclude that $F(x) \leq F(y) .$
\item $F(0)=0$ and $F(1)=1$ follows
quickly from the fact that 0 and 1 are fixed by every $F_{n} .$
\item Finally, if $x$ is a point in $ \mathcal{C}^c,$ then $x$ must fall in the complement of some $\mathcal{C}_{n}.$ Notice that $F(x)=F_{n}(x)$ and the recursive way that each $F_{k}$ is constructed means that in fact $F(y)=F_{k}(y)$ for all $k \geq n$ and $y \in[0,1] \backslash \mathcal{C}_{n}.$
It follows that $F$ is constant on $[0,1]\backslash \mathcal{C}_{n}.$ Thus, $F^{\prime}(x)=0$ for all $x$ in the open set $[0,1]\backslash \mathcal{C}$. Since the Cantor set $\mathcal{C}$ has measure zero $0,$ somehow, $F$ manages to increase from $0$ to $1$ while remaining constant on a set of measure $1$. Thus as we have already mentioned the Cantor-Lebesgue function is an example of a non-constant, {\em singular function} on [0,1].
\end{itemize}

\section{The Lebesgue  space filling curves and beyond}

We have seen that the Cantor-Lebesgue function $F:  \mathcal{C} \to [0,1]$ is a continuous function and it can be extended continuously to $[0,1].$ It is of interest to us that there are continuous mappings from $\mathcal{C}$ onto $[0,1]^2$ and onto $[0,1]^3,$ defined as follows
\begin{eqnarray}\label{LebesgueF2}
F^2\left(0._{3}\varepsilon_{1}\varepsilon_{2}\varepsilon_{3} \ldots \right)=\left(\begin{array}{l}
0._{2}\left(\varepsilon_{1}/2\right)\left(\varepsilon_{3}/2\right)\left(\varepsilon_{5}/2\right) \ldots \\
0._{2}\left(\varepsilon_{2}/2\right)\left(\varepsilon_{4}/2\right)\left(\varepsilon_{6}/2\right) \ldots
\end{array}\right)
\end{eqnarray}
and
\begin{eqnarray}\label{LebesgueF3}
F^3\left(0._{3}\varepsilon_{1}\varepsilon_{2}\varepsilon_{3} \ldots \right)=\left(\begin{array}{l}
0._{2}\left(\varepsilon_{1}/2\right)\left(\varepsilon_{4}/2\right)\left(\varepsilon_{7}/2\right)  \\
0._{2}\left(\varepsilon_{2}/2\right)\left(\varepsilon_{5}/2\right)\left(\varepsilon_{8}/2\right)  \ldots \\
0._{2}\left(\varepsilon_{3}/2\right)\left(\varepsilon_{6}/2\right)\left(\varepsilon_{9}/2\right) \ldots
\end{array}\right)
\end{eqnarray}
where $\varepsilon_{j}=0$ or $2 .$ 

\begin{itemize}
\item Again, as in the case of $F$,  the mappings $F^2: \mathcal{C} \rightarrow [0,1]^2$ and $F^3:\mathcal{C} \rightarrow [0,1]^3$ are not injective but they are surjective. We will prove this for $F^2$ only. The proof for $F^3$ is completely analogous. 
Let $p \in [0,1]^2$. Then, the coordinates of $p$ may be represented by the binary numbers
\[
p=\left(\begin{array}{l}
0._{2} a_{1} a_{2} a_{3} \ldots \\
0._{2} b_{1} b_{2} b_{3} \ldots
\end{array}\right)
\]
If we let $x=0._{\dot{3}}\left(2 a_{1}\right)\left(2 b_{1}\right)\left(2 a_{2}\right)\left(2 b_{2}\right)\left(2 a_{3}\right)\left(2 b_{3}\right) \ldots,$ then by definition  $F^2(x)=p.$
\item The mappings $F^2: \mathcal{C} \rightarrow [0,1]^2$ and $F^3: \mathcal{C} \rightarrow [0,1]^3$ are continuous. Again, we only give a proof for $F^2.$ As done above for $F$, suppose that $\left|x-x_{0}\right|<\frac{1}{3^{2 n}},$ then, as proved above $x_{0}$ and $x$ cannot differ in the first $2 n$ ternary places and we have
\[
\begin{array}{c}
x_{0}=0._{3}\varepsilon_{1}\varepsilon_{2} \ldots\varepsilon_{2 n-1}\varepsilon_{2 n}\varepsilon_{2 n+1} \ldots \\
x=0._{3}\varepsilon_{1}\varepsilon_{2}\ldots\varepsilon_{2 n-1}\tau_{2 n}\tau_{2 n+1} \ldots
\end{array}
\]
\item The mappings $F^2: \mathcal{C} \rightarrow [0,1]^2$ and $F^3: \mathcal{C} \rightarrow [0,1]^3$ are nowhere differentiable. Again, the proof will only be carried out for $F^2.$ As before,
et $x, x_n \in \mathcal{C}$ such that
\[
\begin{aligned}
x &=0._3\varepsilon_{1}\varepsilon_{2} \ldots\varepsilon_{2 n-1}\varepsilon_{2 n}\varepsilon_{2 n+1} \varepsilon_{2 n+2}\ldots \\
x_{n} &=0._3\varepsilon_{1}\varepsilon_{2} \ldots\varepsilon_{2 n-1}\varepsilon_{2 n}\tau_{2 n+1} \varepsilon_{2 n+2} \ldots
\end{aligned}
\]
where $\tau_{2 n+1}=\varepsilon_{2 n+1}+1(\bmod 2) .$ Then
\[
\left|x-x_{n}\right|=2 / 3^{2 n+1}
\]
If $\varphi$ denotes the first component of $F^2$ we have
\[
\varphi(t)-\varphi\left(t_{n}\right)=\left(t_{2 n+1}-\tau_{2 n+1}\right) / 2^{n+1}
\]
and hence,
\[
\left|\varphi(t)-\varphi\left(t_{n}\right)\right| /\left|t-t_{n}\right|=\frac{3}{4}\left(\frac{9}{2}\right)^{n} \longrightarrow \infty,
\]
 thus, $\varphi$ is nowhere differentiable on $\mathcal{C}$. An analogous proof applies to the second component $\psi$ of $F^2.$
\end{itemize}

In $1904,$ H. Lebesgue extended the mapping $F^2$ continuously into $[0,1]$  by linear interpolation. If $\left(a_{n}, b_{n}\right)$ is an interval that is removed in the construction of $\mathcal{C}$ at the $n$-th step, then the extended function $F^2_{l}$ is defined in terms of $F^2$ on that interval as follows:
\begin{equation}\label{LebesgueFL}
F^2_{l}(x)=\frac{1}{b_{n}-a_{n}}\left[F^2\left(a_{n}\right)\left(b_{n}-x\right)+F^2\left(b_{n}\right)\left(x-a_{n}\right)\right], a_{n} \leq x \leq b_{n}
\end{equation}
By construction, $F^2_{l}$ is continuous on $\mathcal{C}^{c}$ and maps $[0,1]$ onto $[0,1]^2$. Moreover, we have 
 \begin{theorem}(Lebesgue) The mapping $F^2_{l}$  defined in $\mathcal{C}$ as  (\ref{LebesgueF2}) and  extended into $\mathcal{C}^{c}$  as (\ref{LebesgueFL}) is continuous on $[0,1],$
  i.e., for any $x_0 \in [0,1]$ we need to prove that given $\varepsilon >0$  there exist $\delta >0$  such  that if  $x \in [0,1]$ such that $|x-x_0|<\delta$ then  
  $$\|F^2_l(x)-F^2_l(x_0)\| \leq \varepsilon .$$
  \end{theorem}
  \dem
  The argument is similar as the one done above for $F$.
\begin{itemize}
\item If $x_{0} \in \mathcal{C}^c,$ then $x_0 \in (a_n, b_n)$ where the open interval $(a_n, b_n)$ is one of those that has been removed from $[0,1]$ in the construction of the Cantor set and then $F^2_{l}$ is trivially at $x_0$ by construction.
\item If $x_{0} \in \mathcal{C},$ then, as $\mathcal{C}$ is perfect, it is either a left accumulation point, a right accumulation point or an accumulation point. We will only consider the case where $x_{0}$ is a left accumulation point, which makes it a right endpoint of an interval that has been removed in the construction of $\mathcal{C}.$ From above we know that namely $x_0 = 0._{3}\varepsilon_{1}\varepsilon_{2}\varepsilon_{3} \ldots \varepsilon_{n-1}2.$ Then we know, that $F^2$ is  the restriction of $F^2_{l}$ to $\mathcal{C},$ is continuous on it and, since $\mathcal{C}$ is compact, uniformly continuous on $\mathcal{C},$ i.e., for any $\epsilon>0,$ there is $\delta>0$ such that
\[
\|F^2\left(x_{1}\right)-F^2\left(x_{2}\right)\|<\epsilon,
\]
for all $x_{1}, x_{2} \in \mathcal{C}$ for which $\left|x_{1}-x_{2}\right|<\delta.$ \\

Let us now examine $|F^2\left(x\right)-F^2\left(x_{0}\right)|.$ There are two possibilities: 
\begin{itemize}
\item $x \in \mathcal{C}$: then, we already know
\[
\|F^2(x)-F^2\left(x_{0}\right)\|<\epsilon \quad \text { for } \quad\left|x-x_{0}\right|<\delta,
\]

\item  $x \in \mathcal{C}^c,$ then, $x \in(a, b),$ where $(a, b)$ is one of the intervals that have been removed in the construction of $\mathcal{C},$ then  $a, b \in \mathcal{C},\;(a, b) \subset \mathcal{C}^{c} .$ But then,
\[
F^2_{l}(x)-F^2_{l}\left(x_{0}\right)=\frac{1}{b-a}\left[\left(F^2(b)-F^2\left(x_{0}\right)\right)(x-a)+\left(F^2(a)-F^2\left(x_{0}\right)\right)(b-x)\right]
\]
 for all $x \in(a, b) .$ Now, let $(a, b)$ be any such interval with $x_{0}<a<b<x_{0}+\delta,$ 
 \[
\|F^2_l(x)-F^2_l\left(x_{0}\right)\|<\frac{1}{b-a}(t-a+b-t) \epsilon=\epsilon
\]
for all $t \in(a, b) \subset \mathcal{C}^{c},$ as  $a \in \mathcal{C}$. \ep
\end{itemize}
\end{itemize}
This, together with the fact that to the left of $x_{0} (x_{0}$ being the right endpoint of a removed interval) $F^2_l$ is continuous by construction, implies the continuity of $F^2_l$ at $x_{0} \in \mathcal{C}.$ An analogous proof applies to a right accumulation point of $\mathcal{C},$ and both proofs together take care of a two-sided accumulation point of $\mathcal{C} .$ \\

In $1878,$ George Cantor demonstrated that any two finite-dimensional smooth manifolds, no matter what their dimensions, have the same cardinality, i.e., it is possible to establish a bijection between them. Cantor's finding implies, in particular, that the interval [0,1] can be mapped bijectively onto the square $[0,1]^{2} .$ The question arose almost immediately whether or not such a mapping can possibly be continuous. In $1879,$ E. Netto showed following theorem. 

\begin{theorem}(Netto): If $f$ represents a bijective map from an $m-$dimensional smooth manifold $\mu_m$ onto an $n-$dimensional smooth manifold $\mu_n$ and $m\neq n$, then $f$ is necessarily discontinuous.
 \end{theorem}
 
 We will prove this theorem for the case where $\mu_m=[0,1]$ and $\mu_n=[0,1]^2$ or $[0,1]^3$. For the general case, we refer the reader to the literature \cite{Green}.
 
\dem  Suppose $f$ is a bijective map from $[0,1]$ to $[0,1]^2$ (or $[0,1]^3$), its inverse $f^{-1}:[0,1]^2\to[0,1]$ exists. Let $g:=f^{-1}$. Let $\mathcal{A}\subseteq\mathbb{R}$ denote a closed set of real line. Then, $\mathcal{A}\cap[0,1]$ is bounded and closed and, hence, compact. Thus, $f(\mathcal{A}\cap[0,1])$,image of $(\mathcal{A}\cap[0,1])$ mapped by $f$, is compact and, hence, closed. Let $\mathcal{A}_1:=f(\mathcal{A}\cap[0,1])$. Since $f(\mathcal{A}\cap[0,1])=g^{-1}(\mathcal{A})$, we have $g^{-1}(\mathcal{A})=\mathcal{A}_1$ and since $\mathcal{A}_1\subseteq f\big([0,1]\big)=\mathcal{D}(g)$ where $\mathcal{D}(g)$ denotes the domain of $g$, we can write $g(\mathcal{A})=\mathcal{A}_1\cap\mathcal{D}(g)$ and we have that $g=f^{-1}$ is continuous.

Since $f^{-1}$ is continuous, it maps connected sets onto connected sets. We remove a point $t_0$ from the open interval $(0,1)$ and its image $f(t_0)$ from $[0,1]^2$ (or $[0,1]^3$). $[0,1]^2\backslash \{f(t_0)\}$ (or $[0,1]^3\backslash \{f(t_0)\}$) are still connected but $[0,1]\backslash\{t_0\}=[0,t_0)\cup(t_0,1]$ is not, and the continuous function $f^{-1}$ appears to map a connected set onto a disconnected set. This is obviously a contradiction and it follows that $f$ cannot be continuous.
\ep\\

Now, suppose the condition of  being bijective is dropped; is it still possible to obtain a continuous surjective mapping from $[0,1]$ onto $[0,1]^{2}$? Since a continuous mapping from $[0,1]$ (or any other interval $(a,b)$) into the plane (or space) was and, to a large extent, still is called a curve, the question may be rephrased as follows: Is there a curve that passes through every point of a two-dimensional region (such as, for example, $[0,1]^{2}$ ) with positive Lebesgue measure? In 1890 G. Peano constructed the first such curve. Curves with this property are now called {\em space-filling curves} or {\em Peano curves.} Further examples  were give by D. Hilbert (in 1891), E. H. Moore (in 1900), W.L. Osgood  (in 1903), H. Lebesgue (in 1904), W. Sierpinski (in 1912), G. P\'olya (in 1913), I.J. Schoenberg (in 1938) and many others. The functions $F^2$ and $F^3$ are called {\em Lebesgue's space filling curve}. For more on space filling curve we refer to \cite{sagan} where all these curves are  carefully studied. \\

Moreover, in 1927, in the second edition of Felix Hausdorff's book {\em Mengenlehre}, there is a amazing result that generalizes the results what we have discussed in this section: 
\begin{theorem}(Hausdorff): Every compact set in $ \mathbb{R}^{n}$ is a continuous image of the Cantor set $\mathcal{C}.$
 \end{theorem}
\dem  Let $\mathcal{K}$ denote a compact set  in $ \mathbb{R}^{n}$. We cover $\mathcal{K}$ by $N_{1}$-neighborhoods of radius $1$. Since $\mathcal{K}$ is compact, there exists a finite subcover:
\[
N_{1}\left(x_{i_{1}}\right), i_{1}=0,1,2, \ldots, 2^{n_{1}}-1
\]
(where we may assume, without loss of generality, that this subcovering has $2^{n_{1}}$ members for some integer $n_{1}-$ by counting some, if necessary, several times). Let
\[
\mathcal{K}_{i_{1}}=\overline{N}_{1}\left(x_{i_{1}}\right) \cap \mathcal{K}
\]
where $\overline{N}_{1}$ denotes the closure of $N_{1},$
and we have $\bigcup_{i_{1}} \mathcal{K}_{i_{1}}=\mathcal{K} .$ 

Next, we cover each $\mathcal{K}_{i_{1}}$ by $N_{1 / 2}$-neighborhoods  of radius $1/2$ and pick finite subcovering $N_{1 / 2}\left(x_{i_{1} i_{2}}\right), i_{2}=0,1,2, \ldots, 2^{n_{2}}-1$ (where we assume, without loss of generality, that each subcovering has the same number of $2^{n_{2}}$ members). Let
$$
\mathcal{K}_{i_{1} i_{2}}=\overline{N}_{1 / 2}\left(x_{i_{1} i_{2}}\right) \cap \mathcal{K}_{i_{1}}
$$

We have $\bigcup_{i_{2}} \mathcal{K}_{i_{1} i_{2}}=\mathcal{K}_{i_{1}} .$ We continue in this manner with $N_{1 / 4}$-neighborhoods of radius $1/4$, $N_{1 / 8}$ -neighborhoods of radius $1/8$,  and so on, obtaining nested sequences of compact sets
\[
\mathcal{K} \supseteq \mathcal{K}_{i_{1}} \supseteq \mathcal{K}_{i_{1} i_{2}} \supseteq \mathcal{K}_{i_{1} i_{2} i_{3}} \supseteq \ldots
\]
where $\mathcal{K}_{i_{1} i_{2} i_{3} \ldots i_{k}}$ is contained in the closure of a $1 / 2^{k-1}$-neighborhood. We may write any $x \in \mathcal{C}$ as
\[
x=0._{3}\varepsilon_{1}\varepsilon_{2} \ldots\varepsilon_{n_{1}}\varepsilon_{n_{1}+1}\varepsilon_{n_{1}+2} \ldots \varepsilon_{n_{1}+n_{2}}  \varepsilon_{n_{1}+n_{2}+1} 
 \ldots
\]
(where as usual $\varepsilon_{j}=0$ or $2$ and where infinitely many trailing $\varepsilon_{j}$ might well be 0) and define a mapping $f: \mathcal{C} \rightarrow \mathcal{K}$ by
\[
f(x)=\mathcal{K}_{i_{1}} \cap \mathcal{K}_{i_{1} i_{2}} \cap \mathcal{K}_{i_{1} i_{2} i_{3}} \cap \mathcal{K}_{i_{1} i_{2} i_{3} i_{4}} \cap \ldots
\]
where 
\begin{eqnarray*}
i_{1}&=&\left(\left( \varepsilon_{n_1}/2\right) \cdots \left(\varepsilon_{n_{1}}/2\right)\right)_{2}, i_{2}=\left(\left(\varepsilon_{n_{1}+1}/2\right) \ldots \left(\varepsilon_{n_{1}+n_{2}}/2\right)\right)_{2},\\
 i_{3}&=&\left(\left(\varepsilon_{n_{1}+n_{2}+1}/2\right) \ldots\left(\varepsilon_{n_{1}+n_{2}+n_{3}}/2\right)\right)_{2}, \ldots
\end{eqnarray*}
with 
$$\left(\left(\varepsilon_{n_{1}+\cdots+n_{k}+1} /2\right)\cdots \left(\varepsilon_{n_{1}+\cdots+n_{k+1}}/2\right)\right)_{2}$$ denotes the binary representation of the number $i_{k+1} .$ By Cantor's intersection theorem, associates with each $x \in \mathcal{C}$ a unique point in $\mathcal{K}$. This mapping is surjective because each point in $\mathcal{K}$ lies in at least one nested sequence which corresponds to a binary sequence and which, in turn, corresponds to a point in $\mathcal{C}$. The mapping is also continuous: If $\left|x^{\prime}-x^{\prime \prime}\right|<1 / 3^{n+1},$ then
\[
\begin{aligned}
x^{\prime}&=0._3\varepsilon_{1}\varepsilon_{2} \ldots\varepsilon_{2 n-1}\varepsilon_{2 n}\varepsilon_{2 n+1} \varepsilon_{2 n+2}\ldots \\
x^{\prime \prime}&=0._3\varepsilon_{1}\varepsilon_{2} \ldots\varepsilon_{2 n-1}\varepsilon_{2 n}\tau_{2 n+1} \tau_{2 n+2} \ldots
\end{aligned}
\]
If $n_{1}+n_{2}+\cdots+n_{j} \leq n<n_{1}+n_{2}+\cdots+n_{j}+n_{j+1},$ then $f\left(x^{\prime}\right), f\left(x^{\prime \prime}\right) \in$
$\mathcal{K}_{i_{1} i_{2} \ldots i_{j}} \subseteq \bar{N}_{1 / 2^{j-1}}\left(x_{i_{1} \ldots i_{j}}\right),$ i.e., 
$$\|f\left(x^{\prime}\right)-f\left(x^{\prime \prime}\right)\|<1 / 2^{j-2}.$$
\ep\\

 Then, from the discussion above there is a natural question: given a compact set $\mathcal{K}$, could the mapping $f: \mathcal{C} \to \mathcal{K},$ obtained by Hausdorff, be extended continuously into the entire interval $[0,1]$ by joining the images of the beginning points and endpoints of the open intervals that make up $\mathcal{C}^{c}?$ This was possible for  the case of $[0,1]^2$ and $[0,1]^3$  because the target set was convex. In general, we need is some guarantee that these points can still be connected by a continuous path that lies in $\mathcal{K},$ even when $\mathcal{K}$ is not convex. H. Hahn and St. Mazurkiewicz, proved that if  that if $\mathcal{K}$ is compact, connected, and locally connected, such a construction is possible. For more details we refer to \cite[Chapter 6]{sagan}.

\end{document}